\documentclass[11pt,a4paper,reqno]{amsart}
\usepackage{amsfonts}
\usepackage{amsmath,amssymb,amsthm,amsxtra}
\usepackage{float}
\usepackage[colorlinks, linkcolor=blue,anchorcolor=Periwinkle,
    citecolor=red,urlcolor=Emerald]{hyperref}
\usepackage[usenames,dvipsnames]{xcolor}
\usepackage{enumitem}
\usepackage{geometry,array} 
\usepackage{graphicx}
\usepackage{subfigure}
\usepackage{bookmark}
\usepackage{tikz}\usetikzlibrary{matrix}
\usepackage{url}

\usetikzlibrary{decorations.markings}
\tikzset{->-/.style={decoration={  markings,  mark=at position #1 with
    {\arrow{>}}},postaction={decorate}}}
\tikzset{-<-/.style={decoration={  markings,  mark=at position #1 with
    {\arrow{<}}},postaction={decorate}}}
\usepackage[all]{xy}
\usetikzlibrary{arrows}
\usetikzlibrary{decorations.pathreplacing,decorations.pathmorphing,shadings,fadings,calc}

\usepackage{extarrows}
\def\dexc{cyan!21}

\theoremstyle{plain}
\newtheorem{theorem}{Theorem}[section]

\newtheorem{lemma}[theorem]{Lemma}

\newtheorem{proposition}[theorem]{Proposition}

\theoremstyle{definition}
\newtheorem{definition}[theorem]{Definition}

\newtheorem{example}[theorem]{Example}

\newtheorem{remark}[theorem]{Remark}

\numberwithin{equation}{section}

\def\hua{\mathcal}
\def\kong{\mathbb}
\def\<{\langle}
\def\>{\rangle}

\def\ZZ{\mathbb{Z}}

\def\RR{\kong{R}}
\def\CC{\kong{C}}

\def\Aut{\operatorname{Aut}}
\def\Ind{\operatorname{Ind}}
\def\Sim{\operatorname{Sim}}
\def\Hom{\operatorname{Hom}}

\def\Stab{\operatorname{Stab}}
\def\Stap{\operatorname{Stab}^\circ}
\def\diff{\operatorname{d}}

\newcommand{\h}{\operatorname{\hua{H}}}            

\renewcommand{\k}{\mathbf{k}}
\renewcommand{\mod}{\operatorname{mod}}
\newcommand{\Ho}[1]{\operatorname{\bf H}_{#1}}
\newcommand{\tilt}[3]{{#1}^{#2}_{#3}}
\newcommand{\Cone}{\operatorname{Cone}}

\def\numbers{\begin{enumerate}[label=\arabic*{$^\circ$}.]}
\def\ends{\end{enumerate}}

\newcommand{\EG}{\operatorname{EG}}       
\newcommand{\C}{\hua{C}}

\newcommand{\D}{\operatorname{\hua{D}}}

\newcommand{\per}{\operatorname{per}}
\newcommand\Sph{\operatorname{Sph}}

\newcommand{\Tri}{\bigtriangleup}

\def\arrow{red}
\def\surf{\mathbf{S}}                       

\newcommand{\ST}{\operatorname{ST}}        
\newcommand{\BT}{\operatorname{BT}}        

\newcommand{\MCG}{\operatorname{MCG}}

\newcommand{\Int}{\operatorname{Int}}

\def\Coh{\operatorname{Coh}}

\def\Fuk{\operatorname{Fuk}}
\newcommand{\DCoh}[1]{\D^b(\Coh #1)}

\def\TT{\kong{T}}
\def\T{T}

\def\M{\mathbf{M}}
\def\P{\mathbf{P}}

\def\surfo{{\mathbf{S}}_\Tri}

\newcommand\Bt[1]{\operatorname{B}_{#1}}

\newcommand{\Quad}{\operatorname{Quad}}

\newcommand\Br{\operatorname{Br}}

\newcommand{\twi}{\Psi} 
\title[Top. stab. and top. Fukaya cat.]
{Topological structure of spaces of stability conditions and
topological Fukaya type categories}

\author{Yu Qiu}
\address{Department of Mathematics,
The Chinese University of Hong Kong, Shatin, Hong Kong.}
\email{Yu.Qiu@Bath.edu}

\date{\today}
\begin{document}

\begin{abstract}
This is a survey on two closely related subjects.
First, we review the study of topological structure of
`finite type' components of spaces of Bridgeland's stability conditions
on triangulated categories \cite{W,KQ,Q2,QW,BPP}.
The key is to understand Happel-Reiten-Smal{\o} tilting as tiling of cells.
Second, we review topological realizations of
various Fukaya type categories \cite{BZ,QZ,QQ,QQ2,QZ2,HKK},
namely cluster/Calabi-Yau and derived categories from surfaces.
The corresponding spaces of stability conditions
are of `tame' nature
and can be realized as moduli spaces of quadratic differentials due to
Bridgeland-Smith and Haiden-Katzarkov-Kontsevich \cite{BS,HKK,I,KQ2}.

\end{abstract}

\keywords{stability conditions, spherical twists, Calabi-Yau categories, quadratic differentials}

\maketitle


\def\dgen{2\mathbb{N}_{\leq g}-1}
\def\dgene{2\mathbb{N}_{\leq g}}

\section{Introduction}
The notion of a stability condition on a triangulated category $\D$
was introduced by Bridgeland \cite{B1}.
The idea was inspired from physics by studying D-branes in string theory.
Bridgeland showed a key result that
the space $\Stab(\D)$ of stability conditions on $\D$
is a finite dimensional complex manifold,
provided that the rank of the Grothendick group $K(\D)$ is finite.
Moreover, these spaces carry interesting geometric/topological structures
which shed light on the properties of the original triangulated categories.

\begin{figure}[t]\centering
\def\graphc{JungleGreen}
\def\vertexc{blue}
\def\halfc{blue}
\def\thirdc{blue}
\def\dexc{black!10!blue!50!green}
\begin{tikzpicture}[scale=.7]
\path (0,0) coordinate (O);
\path (0,6) coordinate (S1);
\draw[fill=\graphc!7,dotted]
    (S1) arc(360/3-360/3+180:360-360/3:10.3923cm) -- (O) -- cycle;
\draw[,dotted] (210:6cm) -- (O);
\draw[fill=black] (30:1.6077cm) circle (.05cm);
\draw
    (O)[Periwinkle,dotted,thick]   circle (6cm);
\draw
    (S1)[\graphc,thick]
    \foreach \j in {1,...,3}{arc(360/3-\j*360/3+180:360-\j*360/3:10.3923cm)} -- cycle;
\draw
    (S1)[\graphc,semithick]
    \foreach \j in {1,...,6}{arc(360/6-\j*360/6+180:360-\j*360/6:3.4641cm)} -- cycle;
\draw
    (S1)[\graphc]
    \foreach \j in {1,...,12}{arc(360/12-\j*360/12+180:360-\j*360/12:1.6077cm)}
        -- cycle;
\foreach \j in {1,...,3}
{\path (-90+120*\j:.7cm) node[\vertexc] (v\j) {$\bullet$};
 \path (-210+120*\j:.7cm) node[\vertexc] (w\j) {$\bullet$};
 \path (-90+120*\j:2.2cm) node[\vertexc] (a\j) {$\bullet$};
 \path (-90+15+120*\j:3cm) node[\vertexc] (b\j) {$\bullet$};
 \path (-90-15+120*\j:3cm) node[\vertexc] (c\j) {$\bullet$};}
\foreach \j in {1,...,3}
{\path[->,>=latex] (v\j) edge[\thirdc,bend left,thick] (w\j);
 \path[->,>=latex] (a\j) edge[\thirdc,bend left,thick] (b\j);
 \path[->,>=latex] (b\j) edge[\thirdc,bend left,thick] (c\j);
 \path[->,>=latex] (c\j) edge[\thirdc,bend left,thick] (a\j);}
\foreach \j in {1,...,3}
{\path (60*\j*2-1-60:3.9cm) node[\vertexc] (x1\j) {$\bullet$};
 \path (60*\j*2+1-120:3.9cm) node[\vertexc] (x2\j) {$\bullet$};
}
\foreach \j in {1,...,3}
{\path[->,>=latex] (v\j) edge[\halfc,bend left,thick] (a\j);
 \path[->,>=latex] (a\j) edge[\halfc,bend left,thick] (v\j);
 \path[->,>=latex] (b\j) edge[\halfc,bend left,thick] (x1\j);
 \path[->,>=latex] (c\j) edge[\halfc,bend left,thick] (x2\j);
 \path[->,>=latex] (x1\j) edge[\halfc,bend left,thick] (b\j);
 \path[->,>=latex] (x2\j) edge[\halfc,bend left,thick] (c\j);}
\end{tikzpicture}
\caption{Exchange graphs as the skeleton of spaces of stability conditions}\label{fig:Q}
\end{figure}
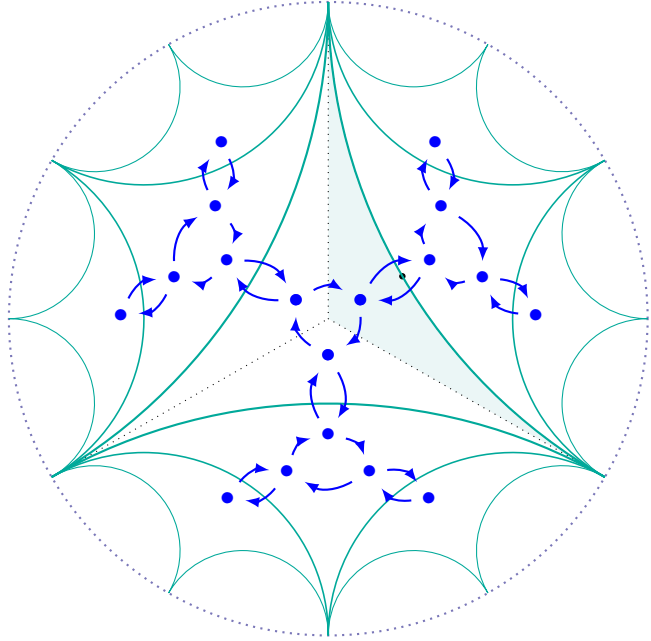

While the spaces of stability conditions in geometric setting are hard to studied in general,
the algebraic (quivery) case is much easier.
Namely, when the heart $\h$ of a stability conditions is finite
(i.e. a length category with finitely many simples),
the stability conditions supported on $\h$ is a complex half open half closed cell,
cf. \eqref{eq:H^n}.
This allows us to study the global topological structure of $\Stab(\D)$.
The skeleton of a finite type component $\Stab_0$
of a space of stability conditions is a (connected component of) exchange graph $\EG_0$.
The vertices and edges of such a graph are (finite) hearts and simple HRS-tilting.
In fact, $\EG_0$ is the dual of some of the stratification of $\Stab_0$
(see \cite[Theorem~3.4]{Q2}, \cite[\S~4.7]{KQ2} and Figure~\ref{fig:Q}).
Applying this method, we show that such a component is always contractible
(Theorem~\ref{thm:QW}).

Next class of components to be studied is the `tame type'.
Our philosophy/definition is that in a tame type component,
any stability condition, up to a $\kong{C}$-action, admits a finite heart.
The corresponding categories, which have been studied, are from surfaces.
The key here is to understand (directly or via other theory, i.e. cluster theory)
objects/morphisms in this type of categories via arcs/intersections on surfaces.
In some sense, they are of topological Fukaya type (cf. Remark~\ref{rem}).

Moreover, these stability spaces
are related to Kontsevich's homological mirror symmetry for a pair of mirror varieties $(X,X^\vee)$.
More precisely,
the K\"{a}hler moduli space $\hua{M}(X^\vee)$ of $X^\vee$
can be (conjecturally) embedded into
(a quotient of) the space of stability conditions on the categories
\[
        \D\Fuk(X) \cong \DCoh{X^\vee},
\]
where $\D\Fuk(X)$ is the derived Fukaya category of $X$ and
$\DCoh{X^\vee}$ the bounded derived category of coherent sheaves on $X^\vee$:

Topological Fukaya type categories are the test ground for this idea/conjecture.
Two pioneer works have been carried out by Bridgeland-Smith and Haiden-Katzarkov-Kontsevich
in \cite{BS,HKK}.
Our upcoming works \cite{IQ1,IQ2} aim to understand the relations between their results.

\subsection*{Context}\label{c}
In \S~\ref{sec:struc}
we recall couple of structures on abelian/triangulated categories.
In \S~\ref{sec:TT}
we review a series of works \cite{W,KQ,QW}, which use HRS-tilting as tiling (and stratification),
to understand the (global) topological structure of finite type components of
spaces of stability conditions.
In \S~\ref{sec:CC}
we review a series of works \cite{BZ,QZ,QQ,QQ2,QZ2,HKK}, which use (graded) string model,
to understand cluster/derived categories via arcs/intersections on surfaces.
In \S~\ref{sec:Quad}
we collect two key results \cite{BS,HKK} in the direction of the study of stability conditions
via quadratic differentials.

\subsection*{Acknowledgements}
This survey is contributed to
the proceedings of the first annual meeting of ICCM = International Consortium of Chinese Mathematicians.
This work is supported by Hong Kong RGC 14300817
(from Chinese University of Hong Kong).

\section{Structures on triangulated categories}\label{sec:struc}
\subsection{Torsion pairs and stability functions on abelian categories}\label{sec:TP}
\begin{definition}
A \emph{torsion pair} in an abelian category $\h$ is an ordered pair of
full subcategories $\<\hua{T}, \hua{F}\>$, 
such that $\Hom(\hua{T},\hua{F})=0$ and any object $E \in \h$ sits
in a short exact sequence
$0 \to E^{\hua{T}} \to E \to E^{\hua{F}} \to 0$,
for some $E^{\hua{T}} \in \hua{T}$ and $E^{\hua{F}} \in \hua{F}$.
We will call $\hua{T}$ the torsion part and $\hua{F}$ the torsion free part.
We will write $\h=\<\hua{T}, \hua{F}\>$ for a torsion pair of $\h$.
\end{definition}

Here is a well-known lemma (cf. \cite[Lemma~3.6]{KQ}.
\begin{lemma}\label{lem:simpletilt}
Let $S$ be a rigid simple object in a Hom-finite abelian category $\h$.
Then $\h$ admits a torsion pair $\<\hua{T},\hua{F}\>$
such that $\hua{F}=\<S\>$.
Similarly, there is also a torsion pair with the torsion part $\hua{T}=\<S\>$.
\end{lemma}

Next we have the notion of stability function on abelian categories,
which is closely related to the corresponding notion on triangulated categories.
Denote by $K(\h)$ the Grothendick group of $\h$.

\begin{definition}
A \emph{stability function} on an abelian category $\h$
is a group homomorphism $Z: K(\h)\to\kong{C}$ such that
for any non-zero object $M$, we have
$Z(M) = m(M) \exp(\phi_Z(M) \mathbf{i} \pi)$ for some
$m(M) \in \kong{R}_{>0}$ and $\phi_Z(M)\in (0,1]$,
i.e. $Z(M)$ lies in the upper half-plane
\begin{gather}\label{eq:H}
    \kong{H}=\{r e^{ \theta  \pi \mathbf{i} }\mid r\in\kong{R}_{>0},0< \theta\le1 \}
    \subset\kong{C}.
\end{gather}
Call $\phi_Z(M)$ the phase of $M$.
We say an object $0\neq M\in\hua{C}$ is semistable (with respect to $Z$) if every
subobject $0\neq L$ of $M$ satisfies $\phi_Z(L) \leq \phi_Z(M)$.
Further, we say a stability function $Z$ on $\hua{C}$ satisfies HN-property,
if for an object $0\neq M\in\hua{C}$,
there is a collection of short exact sequences
\begin{equation*}
\xymatrix@C=1pc@R=1.4pc{
  0=M_0 \ar@{^{(}->}[rr] && M_1 \ar@{->>}[dl] \ar@{^{(}->}[rr] &&  ...
  \ar@{^{(}->}[rr] && M_{k-1} \ar@{^{(}->}[rr] && M_k=M \ar@{->>}[dl] \\
  & L_1 && && && L_k
}\end{equation*}
in $\hua{C}$ such that $L_1,...,L_k$ are semistable objects
(with respect to $Z$) and their phases are in decreasing order,
i.e. $\phi(L_1)>\cdots>\phi(L_k)$.
\end{definition}

The stability functions on abelian categories can be used
to define (motivic) Donaldson-Thomas invariants via quantum dilogarithm identities
(due to Reineke and Kontsevich-Soibelman, cf. the survey \cite{K6}).
One can also use Keller's green mutation to define such invariants combinatorially.
\cite[\S~8]{Q2} provides an alternative proof of quantum dilogarithm identities
in the Dynkin case via topological structure of spaces of stability conditions
that will be discussed in the following.

\subsection{Stability structures on triangulated categories}\label{sec:t-str}
A \emph{torsion pair} in a triangulated category $\D$ is an ordered pair of
full subcategories $\<\hua{T}, \hua{F}\>$, 
such that $\Hom(\hua{T},\hua{F})=0$ and any object $E \in \D$ sits
in a triangle
$$ E^{\hua{T}} \to E \to E^{\hua{F}} \to E^{\hua{T}}[1],$$
for some $E^{\hua{T}} \in \hua{T}$ and $E^{\hua{F}} \in \hua{F}$.
As triangulated categories admit distinguish auto-equivalence, the shift $[1]$,
we have the following upgraded version of torsion pairs.
\begin{definition}
A t-structure $\hua{P}$ is the torsion part of a torsion pair satisfying
$\hua{P}[1]\subset\hua{P}$.

The \emph{heart} of a t-structure $\hua{P}$ is the full subcategory
\[  \h=  \hua{P}^\perp[1]\cap\hua{P}, \]
where $\hua{P}^\perp$ is the corresponding torsion free part.
\end{definition}

A t-structure $\hua{P}$ is \emph{bounded} if
\[
  \hua{D}= \displaystyle\bigcup_{i,j \in \ZZ} \hua{P}^\perp[i] \cap \hua{P}[j].
\]
Equivalently, one can define the heart as follows.
\begin{definition}
A heart $\h$ of a triangulated category $\D$ is
an abelian subcategory satisfying the following conditions:
\begin{itemize}
\item[(1)]if $k_1 > k_2 (\in\ZZ)$ and $A_i \in \h[k_i]\,(i =1,2)$, then $\Hom(A_1,A_2) = 0$,
\item[(2)]for $0 \neq E \in \D$,
there is a finite sequence of integers
\begin{equation*}
k_1 > k_2 > \cdots > k_m
\end{equation*}
and a collection of exact triangles (known as HN-filtration)
\begin{equation}\label{eq:HN}
0 =
\xymatrix @C=5mm{
 E_0 \ar[rr]   &&  E_1 \ar[dl] \ar[rr] && E_2 \ar[dl]
 \ar[r] & \dots  \ar[r] & E_{m-1} \ar[rr] && E_m \ar[dl] \\
& A_1 \ar@{-->}[ul] && A_2 \ar@{-->}[ul] &&&& A_m \ar@{-->}[ul]
}
= E
\end{equation}
with $A_i \in \h[k_i]$ for all $i$.
\end{itemize}
\end{definition}

It is well-known that
\begin{itemize}
    \item there is a bijection between the set of bounded t-structures
    and the set of hearts on a triangulated category, sending a t-structure to its heart.
\end{itemize}
Note that \eqref{eq:HN} means that a heart $\h$ provides
a homology (or a $\ZZ$-structure) on $\D$, that
one can define the $k$-th homology of $E$, with respect to $\h$,
to be
\begin{equation}\label{eq:homology}
 \Ho{k}(E)=
 \begin{cases}
   A_i & \text{if $k=k_i$} \\
   0 & \text{otherwise.}
 \end{cases}
\end{equation}
Then $\hua{P}$ consists of those objects
with no (nonzero) negative homology,
$\hua{P}^\perp$ those with only negative homology
and $\h$ those with homology only in degree 0.

One can refine the (bounded) t-structure ($\ZZ$-structure) on triangulated categories
by further dividing, i.e. to get the $\RR$-structure as follows.

\begin{definition}
A \emph{slicing} on $\D$ is a collection of additive subcategories $$\{\hua{P}(\phi)\subset\D\mid\phi\in\RR\}$$
such that, if one defines (by triangulated extension)
$$\hua{P}(I)=\<\hua{P}(\phi)\mid\phi\in I\>$$
for an interval $I\subset\RR$, then
\begin{itemize}
\item $\hua{P}(\phi+1)=\hua{P}(\phi)[1], \forall \phi\in\RR$.
\item $\hua{P}[\phi,+\infty)$ and $\hua{P}(\phi,+\infty)$ are t-structures for any $\phi\in\RR$.
\end{itemize}
\end{definition}

Here is Bridgeland's $\CC$-upgrading of slicing on triangulated categories.
\begin{definition}\cite{B1}
A \emph{stability condition} $\sigma = (Z,\hua{P})$ on $\hua{D}$ consists of
a group homomorphism $$Z\colon K(\hua{D}) \to \kong{C}$$ called the \emph{central charge} and
a slicing $\{\hua{P}(\phi)\mid\phi\in\RR\}$ satisfying the following compatible condition:
\begin{itemize}
\item if $0 \neq E \in \hua{P}(\phi)$
then $Z(E) = m(E) e^{\phi  \pi \mathbf{i} }$ for some $m(E) \in \kong{R}_{>0}$.
\end{itemize}
The (non-zero) objects in the slicing are said to be semistable
and the simple objects in (the abelian category) $\hua{P}(\phi)$ are said to be stable.
\end{definition}
Note that a slicing does not necessarily admit a compatible central charge to
form a stability conditions.
We will always assume the \emph{support property} on stability conditions \cite{KoSo},
i.e. for some norm $\|\cdot\|$ on $K(\D)\otimes \RR$,
there is a constant $C > 0$ such that
$$\|\alpha\|< C ¡¤ |Z(\alpha)|$$
for all classes $\alpha\in K(\D)$ represented by $\sigma$-semistable objects in $\D$.

\begin{theorem}\label{thm:B1}\cite[Theorem~1.2]{B1}
All stability conditions on a triangulated category
$\D$ form a complex manifold, denoted by $\Stab(\D)$;
each connected component of $\Stab(\D)$ is locally homeomorphic to a linear sub-manifold of
$\Hom_{\kong{Z}}(K(\D),\kong{C})$, sending a stability condition $(Z,\hua{P})$ to its central change $Z$.
\end{theorem}

There is a natural $\CC$-action
on the set $\Stab(\D)$ of all stability conditions on $\D$, namely:
\[
    s \cdot (Z,\hua{P})=(Z \cdot e^{-\mathbf{i} \pi s},\hua{P}_{\operatorname{Re}(s)}),
\]
where $\hua{P}_x(\phi)=\hua{P}(\phi+x)$.
There is also a natural action on $\Stab(\D)$ induced by
the auto-equivalence group $\Aut\D$, namely:
$$\Phi  (Z,\hua{P})=(Z \circ \Phi, \Phi (\hua{P})).$$
The heart $\h(\sigma)$ of a stability condition $\sigma=(Z,\hua{P})$ is $\hua{P}(0,1]$.
And we will say $\sigma$ is supported on $\h$ if $\h(\sigma)=\h$.

\section{Tilting as tiling}\label{sec:TT}
Following \cite{W, KQ, Q2, QW}, we review the idea of
(simple) HRS-tilting as tiling
to study finite type components of spaces of stability conditions.
Let $\D$ be a triangulated category as above with $K(\D)\cong\ZZ^n$.
\subsection{HRS-tilting}\label{sec:HRS}
The canonical partial order on the set of hearts in $\D$
are defined by
\[
    \h_1\leq \h_2 \Longleftrightarrow \hua{P}_1 \supset \hua{P}_2£¬
\]
where $\hua{P}_i$ are the corresponding t-structures.
\begin{theorem}\cite{HRS}
For any torsion pair of a heart $\h=\<\hua{T},\hua{F}\>$ (in $\D$)
there exists the following two hearts with torsion pairs
\[
    \h^\sharp=\<\hua{F}[1],\hua{T}\>,\quad \h^\flat=\<\hua{F},\hua{T}[-1]\>.
\]
We call $\h^\sharp$ the \emph{forward tilt} of $\h$
with respect to the torsion pair $\<\hua{F},\hua{T}\>$,
and $\h^\flat$ the \emph{backward tilt} of $\h$.
\end{theorem}

Clearly $\hua{T}=\h\cap\h^\sharp$ and $\hua{F}=\h\cap\h^\sharp[-1]$.
Also we have $\h^\flat=\h^\sharp[-1]$ and
$$\h[-1] \leq \h^\flat \leq \h \leq \h^\sharp \leq \h[1].$$
It is well-known that there are canonical bijections between
\begin{itemize}
\item the set of torsion pairs in $\h$,
\item the set of hearts in $\D$ between $\h[-1]$ and $\h$,
\item the set of hearts in $\D$ between $\h$ and $\h[1]$,
\end{itemize}
which send a torsion pair to its backward and forward tilts respectively.

\begin{definition}\cite{KQ}
A forward tilt of a heart $\h$ is \emph{simple},
if, in the corresponding torsion pair $\<\hua{T},\hua{F}\>$, the 
part $\hua{F}$ is generated by a single rigid simple $S$
(and thus $\hua{T}=^\perp S$).
We denote the new heart by $\tilt{\h}{\sharp}{S}$.
Similarly, a backward tilt of $\h$ is simple if 
$\hua{T}$ is generated by a rigid simple $S$
(and thus $\hua{F}=S^\perp$.)
The new heart is denoted by $\tilt{\h}{\flat}{S}$.

The \emph{total exchange graph} $\EG(\D)$ of a triangulated category $\D$
is the oriented graph
whose vertices are all hearts in $\D$
and whose edges correspond to simple forward tiltings between them.
\end{definition}

For instance, Lemma~\ref{lem:simpletilt} says when a simple is rigid,
then the corresponding forward/backward simple tilts exist.
One of the most important feature of simple tilting is that
one can iterated tilt a (rigid) simple (up to shifts).
Let $\h$ be a heart in a triangulated category $\D$.
Denote by $\Sim\h$ the set of simples in an abelian category $\h$.
For $S \in \Sim\h$, we set $\tilt{\h}{0\sharp }{S}=\h$ and inductively define
(\cite[Def.~5.11]{KQ})
\[
    \tilt{\h}{m\sharp }{S}
    ={  \Big( \tilt{\h}{ (m-1) \sharp}{S} \Big)  }^{\sharp}_{S[m-1]} \, ,
\]
for $m\geq1$, and similarly $\tilt{\h}{m\flat }{S}$, for $m\geq1$.
For $m<0$, we also set $\tilt{\h}{m\sharp }{S}=\tilt{\h}{-m\flat }{S}$.

The formulae to calculate new simples in a forward/backward simple tilt is given in \cite[Prop.~5.4]{KQ}.

\subsection{Finite type contractible connected components}\label{sec:QW}
\newcommand{\cub}{\operatorname{U}} 
\newcommand{\Cub}[1]{\overline{\cub_{#1}}} 

A finite/algebraic heart $\h$ is a length category with finitely many simples.
In this case, the stability conditions that supported on $\h$,
denoted by $\cub(\h)$ is
\begin{equation}\label{eq:H^n}
    \cub(\h)\cong\kong{H}^n,
\end{equation}
where the coordinates are given by the central charge $\{Z(S)\mid S\in\Sim\h\}$
on the simples of $\h$.

\begin{definition}
A finite component $\Stab_0$ in a space $\Stab(\D)$ of stability conditions
is the union
$$\Stab_0=\bigcup_{\h\in\EG_0}\cub(\h)$$
for some connected component $\EG_0$ of $\EG\D$ consisting of hearts that are all finite.
\end{definition}

We have the following result.
\begin{theorem}\cite{QW}\label{thm:QW}
A finite component $\Stab_0$ is a contractible connected component of $\Stab(\D)$.
\end{theorem}

There are the following three of classes of finite components.
For simplicity, let $\k$ be a fixed algebraically-closed field.
Recall that a triangulated category $\D$ is called Calabi-Yau-$N$ (CY-$N$)
if, for any objects $L,M$ in $\D$ we have a natural isomorphism
\begin{equation}\label{eq:serre}
    \mathfrak{S}:\Hom_{\hua{D}}^{\bullet}(L,M)
        \xrightarrow{\sim}\Hom_{\hua{D}}^{\bullet}(M,L)^\vee[N].
\end{equation}

\begin{example}\label{ex:CY-N}
Let $Q$ be a Dynkin quiver (in general $Q$ can be an acyclic quiver).
Denote by $\Gamma_N Q$ the \emph{Calabi-Yau-$N$ Ginzburg (dg) algebra} associated to $Q$,
which is constructed as follows (\cite[\S~7.2]{K10}, \cite{G}):
\begin{itemize}
\item   Let $Q^N$ be the graded quiver whose vertex set is $Q_0$
and whose arrows are: the arrows in $Q$ with degree $0$;
an arrow $a^*:j\to i$ with degree $2-N$ for each arrow $a:i\to j$ in $Q$;
a loop $e^*:i\to i$ with degree $1-N$ for each vertex $e$ in $Q$.
\item   The underlying graded algebra of $\Gamma_N Q$ is the completion of
the graded path algebra $\k Q^N$ in the category of graded vector spaces
with respect to the ideal generated by the arrows of $Q^N$.
\item   The differential of $\Gamma_N Q$ is the unique continuous linear endomorphism homogeneous
of degree $1$ which satisfies the Leibniz rule and
takes the following values on the arrows of $Q^N$:
\[
    \diff a^*=0,\qquad
    \diff \sum_{e\in Q_0} e^*=\sum_{a\in Q_1} \, [a,a^*] .
\]
\end{itemize}
Our primary example is when $\D=\D_{fd}(\Gamma_N Q)$,
the finite dimensional derived category of $\Gamma_N Q$,
and $\Stab_0=\Stab(\D)$.
\end{example}

\begin{example}
Another two classes of finite type components are (cf. \cite{QW,BPP})
\begin{itemize}
  \item $\D$ is a discrete derived category of finite global dimension
  and $\Stab_0=\Stab(\D)$,
  which includes the cases when $\D=\D^b(Q)$ for a Dynkin quiver $Q$.
  \item $\D$ is a locally-finite triangulated category with finite rank Grothendieck
    group and $\Stab_0$ is any connected component of $\Stab(\D)$.
\end{itemize}
\end{example}

\subsection{Tiling and stratification}\label{sec:tiling}
In this section, we explain the idea of the proof of Theorem~\ref{thm:QW}.
We will fix a finite component $\Stab_0$ that corresponds to $\EG_0$.

First, each finite heart corresponds to a cell in $\Stab(\D)$.
\begin{proposition}\cite{B1}\label{pp:ss}
To give stability condition on a triangulated category $\hua{D}$
is equivalent to giving a heart with
a stability function satisfying the HN-property.
Further, to give a stability condition on $\hua{D}$ with a finite heart $\h$
is equivalent to give a function $\Sim\h \to \kong{H}$,
where $\kong{H}$ is the upper half plane as in \eqref{eq:H}.
\end{proposition}

Moreover, simple HRS-tilting corresponds to tiling (cf. \cite{B1, W, Q2})
in the sense as below.

\begin{lemma}
Let $\h_1,\h_2$ be two finite hearts.
Then $\cub(\h_1)\cap\cub(\h_2)$ is a codimension one wall of the cells $\cub(\h_i)$
if and only if $\h_1$ is a simple tilt of $\h_2$.
\end{lemma}

This implies the following.
Here $\overline{\cub}$ is the closure of $\cub$.

\begin{proposition}\cite[Thm.~3.4]{Q2}\label{pp:Q2}
There is a canonical embedding $\iota:\EG_0 \hookrightarrow \Stab_0$,
unique up to homotopy, such that
\begin{itemize}
\item   $\iota(\h)\in(\cub(\h))^\circ$ for any heart $\h\in\EG_0$.
\item   $\iota(s)$ is contained in $(\cub(\h)\cup\cub(\tilt{\h}{\sharp}{S}))^\circ$
and transversally intersects $$(\cub(\h)\cap\overline{\cub(\tilt{\h}{\sharp}{S})})^\circ$$
at exactly one point,
for any edge $s:\h\to\tilt{\h}{\sharp}{S}$.
\end{itemize}
Further, there is a surjection $\pi_1\EG_0\twoheadrightarrow\pi_1\Stab_0$.
\end{proposition}

Pushing this idea to the limit,
we can construct a stratification of $\Stab_0$ in \cite{QW}.
\newtheorem{construction}[theorem]{Construction}
\begin{construction}
For $\h\in\EG_0$ and $I\subset\Sim\h$, define the stratum
\[
    \cub_{\h,I}\colon=\{ \sigma=(Z,\hua{P})\in\cub(\h) \mid
        \hua{P}(1)=\<I\> \}
        \subset \cub(\h).
\]
Then $\cub_{\h,I}\cong\kong{H}^{n-i}\times\RR_{<0}^i$ for $i=|I|$.
Of course we have
\[
    \Stab_0=\bigcup_{\h\in\EG_0}\cub(\h)
    =\bigcup_{\h\in\EG_0} \left( \bigcup_{I\subset\Sim\h} \cub_{\h,I} \right)
\]
\end{construction}

Moreover,
\[
    \Cub{\h,I}=\bigsqcup_{I\subset K\subset\Sim\h} \partial_K\cub_{\h,I}
\]
for
\[
    \partial_K\cub_{\h,I}\colon=\{ \sigma=(Z,\hua{P})\in\Cub{\h,I} \mid
        Z(S)\in\RR\smallsetminus\{0\} \Longleftrightarrow S\in K \}.
\]
The important observation is the following, which encodes stratification structure of $\Stab_0$.

\begin{proposition}\cite[Cor.~3.13]{QW}
Let $\h,\hua{G}\in\EG_0$ and $I\subset\Sim\h, J\subset\Sim\hua{G}$.
\begin{gather}\label{eq:W}
    \cub_{\h,I}\subset\Cub{\hua{G},J}\;\Longleftrightarrow\;
    \tilt{\h}{\flat}{I}\le \tilt{\hua{G}}{\flat}{J}\le \hua{G}\le \h
\end{gather}
where $\tilt{\h}{\flat}{I}$ is the backward tilt of $\h$ w.r.t.
the torsion pair whose torsion part is generated by $I$
and similarly for $\tilt{\hua{G}}{\flat}{J}$.
\end{proposition}

This implies that the stratum of $\Stab_0$ makes it
a (regular and totally-normal) CW-cellular stratified space,
see \cite{QW} for details.
Then one can prove $\Stab_0$ is contractible as in Theorem~\ref{thm:QW}.

\begin{remark}
Proposition~\ref{pp:Q2} also holds in some general setting,
c.f \cite{BS} and \cite[\S~4.3]{KQ2} in the marked surface case.
\end{remark}

\subsection{On spherical twists: Dynkin case}\label{sec:ST}
One of most crucial classes of auto-equivalences is spherical twists \cite{ST},
which provides braid like actions on $\Stab(\D)$.
\begin{definition}
An object $S$ in a CY-$N$ category $\D$ is \emph{($N$-)spherical} if
$\Hom^{\bullet}(S, S)=\k \oplus \k[-N]$ and induces
a \emph{twist functor} $\phi_S\in\Aut\D$, such that \[
    \twi_S(X)=\Cone\left(S\otimes\Hom^\bullet(S,X)\to X\right)
\]
with inverse
\[
    \twi_S^{-1}(X)=\Cone\left(X\to S\otimes\Hom^\bullet(X,S)^\vee \right)[-1].
\].

In the case when $\D=\D_{fd}(\Gamma_N Q)$,
any simple in the canonical heart $$\h_Q^N\colon=\mod\Gamma_N Q.$$ is spherical and
the \emph{spherical twist group} is
\[
    \ST(\Gamma_N Q)\colon=\< \twi_S \mid S\in\Sim\h_Q^N \>\subset\Aut\D.
\]
We have the following result.
See the survey \cite{Qs} for the full credit on this topic
(generalization of braid group, in particular, spherical twist group).
\end{definition}
\begin{theorem}\cite{QW}
If $Q$ is a Dynkin quiver and $N\ge2$,
then $\ST(\Gamma_N Q)$ is isomorphic to the braid group $\Br(Q)$,
sending the standard generators to the corresponding ones,
and acts faithfully on $\Stab(\D_{fd}(\Gamma_N Q))$.
\end{theorem}

\section{Topological models for categories}\label{sec:CC}

\subsection{Quivers with potential}\label{sec:QP}
Let $(Q,W)$ be a quiver with potential, that is,
a directed graph $Q$ and a linear combination of cycles in $Q$.
\begin{definition}\cite{G, K10}
The \emph{Ginzburg dg algebra (of degree 3)} $\Gamma(Q,W)$ associated to $(Q,W)$
is constructed as follows:
\begin{itemize}
\item $Q^3$ is constructed in Example~\ref{ex:CY-N}
that provides underlying graded algebra structure of $\Gamma(Q,W)$.
\item  The differential of $\Gamma(Q,W)$ is the unique continuous linear endomorphism homogeneous
of degree $1$ which satisfies the Leibniz rule and takes the following values
\[
  \diff \sum_{e\in Q_0} e^*  =  \sum_{a\in Q_1} \, [a,a^*],
  \qquad
  \diff \sum_{a\in Q_1} a^* =  \partial W,
\]
where $\partial W$ is the full cyclic derivative of $W$ with respect to the arrows.
\end{itemize}
\end{definition}

There are three triangulated categories associated to $\Gamma=\Gamma(Q,W)$.
\begin{itemize}
\item the \emph{finite dimensional derived category} $\D_{fd}(\Gamma)$ as mentioned above,
which is CY-3;
\item the \emph{perfect derived category} $\per\Gamma$;
\item the \emph{cluster category} $\C(\Gamma)$ (which is CY-2)
that fits into Amiot's short exact sequence of triangulated categories (\cite{K10})
\begin{gather}
    0 \to \D_{fd}(\Gamma) \to \per\Gamma  \xrightarrow{\pi} \C(\Gamma) \to 0.
\end{gather}
Denote the quotient map by $\pi$.
\end{itemize}

\subsection{Marked surfaces}\label{sec:MS}
\def\CA{\operatorname{CA}}
\def\OA{\operatorname{OA}}
\def\GA{\operatorname{GA}}

Following \cite{FST},
let $\surf$ be a \emph{marked surface} with non-empty boundary
with a finite set $\M$ of marked points on its boundary $\partial \surf$
and a finite set $\P$ of punctures in its interior $\surf\setminus\partial\surf$ such that the following conditions hold
\begin{itemize}
\item each connected component of $\partial \surf$ contains at least one marked point,
\item the rank
\begin{equation}\label{eq:rank}
n=6g+3p+3b+m-6
\end{equation}
of the surface is positive, where $g$ is the genus of $\surf$,
$b$ the number of boundary components,
$m=|\M|$ the number of marked points and $p=|\P|$ the number of punctures.
\end{itemize}
For the most part of the paper, we will assume $\mathbf{P}=\emptyset$
for simplicity. The exceptional case (i.e. results are proved in general setting) will be pointed out.

\begin{definition}
An \emph{open arc} in $\surf$ is a curve (up to isotopy) in its interior,
whose endpoints are marked points in $\M$, and which is neither homotopic to a boundary segment nor to a point.
Denote by $\OA(\surf)$ the set of open arcs on $\surf$.

A triangulation $\T$ of $\surf$
is a maximal collection of open arcs in $\surf$ which do not
intersect each other in the interior of $\surf$.
It is well-known that any triangulation $\T$ of $\surf$ consists of $n$ open arcs and divides $\surf$ into
\begin{gather}\label{eq:aleph}
    \aleph=\frac{2n+|\M|}{3}
\end{gather}
triangles.
\end{definition}

There are the following correspondences between
topological data of $\surf$ and the associated cluster algebras,
which work even when $\mathbf{P}\ne\emptyset$
with certain modification (i.e. tagging).
See \cite{FST} for details.

\begin{table}[h]
\caption{Correspondences I}
\begin{tabular}{ccc}
\\Topology && Cluster Algebra\\[1ex]
\hline
(Tagged) Open arcs in $\surf$ & $\xymatrix@C=2pc{\ar@{<->}[r]^{\text{\tiny{1-1}}}&}$ & Cluster variables\\
\hline
(Tagged) Triangulations of $\surf$ & $\xymatrix@C=2pc{\ar@{<->}[r]^{\text{\tiny{1-1}}}&}$ & Clusters\\
\hline
(Tagged) Flip & $\xymatrix@C=2pc{\ar@{<->}[r]^{\text{\tiny{1-1}}}&}$ & Mutation\\
\hline\hline
\end{tabular}
\end{table}

Let $\T$ be a triangulation of $\surf$.
Then there is an associated quiver $Q_\T$ with a potential $W_\T$, constructed as follows
(see, e.g. \cite{FST} and \cite{LF} for the precise definition):
\begin{itemize}
\item the vertices of $Q_\T$ are (indexed by) the arcs in $\T$;
\item for each triangle in $\T$, there are three arrows between the corresponding vertices as shown in Figure~\ref{fig:quiver};
\item these three arrows form a 3-cycle in $Q_\TT$ and $W_\T$ is the sum of all such 3-cycles.
\end{itemize}
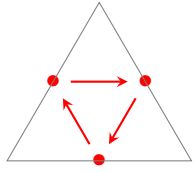
\begin{figure}[h]\centering
  \begin{tikzpicture}[scale=.7]
  \foreach \j in {1,...,3}  { \draw (120*\j-30:2) coordinate (v\j);}
    \path (v1)--(v2) coordinate[pos=0.5] (x3)
              --(v3) coordinate[pos=0.5] (x1)
              --(v1) coordinate[pos=0.5] (x2);
    \foreach \j in {1,...,3}{\draw (x\j) node[red] (x\j){$\bullet$};}
    \draw[->,>=stealth,red,thick] (x1) to (x3);
    \draw[->,>=stealth,red,thick] (x3) to (x2);
    \draw[->,>=stealth,red,thick] (x2) to (x1);
    \draw[gray,thin] (v1)--(v2)--(v3)--cycle;
  \end{tikzpicture}
\caption{The (sub-)quiver (with a potential) associated to a triangle}
\label{fig:quiver}
\end{figure}

\subsection{Cluster categories from marked surfaces}\label{sec:ClusterC}

We have the following result for the cluster categories $\hua{C}(\Gamma_\T)$.
For the unpuncture case, it is due to \cite{BZ} and the general case
($\mathbf{P}\ne\emptyset$) is due to \cite{QZ}.
Note that different triangulations will give the same category up to equivalence.
By abuse of notation, we will write $\hua{C}(\surf)$ for $\hua{C}(\Gamma_\T)$ sometimes.

\begin{theorem}\cite{BZ,ZZZ,QZ}\label{thm:BZZZZQZ}
There is a bijection
\begin{gather}\label{eq:M}
\begin{array}{rcl}
    M\colon\OA(\surf) &\to& \hua{C}^\circ(\surf),\\
    \gamma &\mapsto& M_\gamma,
\end{array}
\end{gather}
where $\hua{C}^\circ(\surf)$ is the set of string objects in the cluster category
$\hua{C}(\surf)$, with the formula
\begin{gather}\label{eq:M int=dim}
    \Int(\gamma_1,\gamma_2)=\dim\Hom^1(M_{\gamma_1},M_{\gamma_2}).
\end{gather}
\end{theorem}
Note that \eqref{eq:M int=dim} is compatible with the fact
that $\hua{C}(\surf)$ is CY-2 as interchanging $\gamma_i$
does not effect the value on both hand sides.
Moreover, the Auslander-Reiten translation $\tau$ on $\hua{C}(\surf)$,
which is the shift functor $[1]$, is realized as
the universal (tagged) rotation (cf. \cite[Def.~3.4]{BQ}) of $\OA(\surf)$.
Then one can apply \eqref{eq:M int=dim} to calculate $\Hom$ of any degree.

\begin{remark}\label{rem}
We consider a category that admits this type of topological realization
is of topological Fukaya type.
More precisely, the topological realization means
\begin{itemize}
  \item An (important) class of objects is bijective to curves/arcs on a surface.
  \item The basis of homomorphism spaces is bijective to intersections between the corresponding curves/arcs.
\end{itemize}
\end{remark}
In the following, we collect couple of results of this type.

\subsection{Derived categories from decorated marked surfaces}\label{sec:DMS}
Recall that any triangulation divides $\surf$ into $\aleph$ triangles,
cf. \eqref{eq:aleph}.
\begin{definition}\cite[Def.~3.1]{QQ}
The \emph{decorated marked surface} $\surfo$ is a marked surface $\surf$
with a set $\Tri$ of $\aleph$ decorating points in
the interior of $\surf$.
Moreover, we have the following notions.
\begin{itemize}
\item A \it{closed arc} in $\surfo$ is a curve (up to isotopy) in
     $\surfo^\circ\smallsetminus\Delta$, whose endpoints are different points in $\Delta$.
     Denote by $\CA(\surfo)$ the set of simple closed curve in $\surfo$.
\item An \emph{open arc} in $\surfo$ is a curve (up to isotopy) in
     $\surfo^\circ\smallsetminus\Delta$, whose endpoints are marked points in $\M$,
     and which is neither homotopic to a boundary segment nor to a point.
    Denote by $\OA(\surfo)$ the set of open arcs on $\surfo$.
\item A triangulation $\TT$ of $\surfo$
is a maximal collection of open arcs in $\surfo$ that divides $\surfo$
into once-decorated triangles.
More precisely, $\surfo$ will be divided into $\aleph$ triangles,
each of which contains exactly one decorating point in $\Delta$.
\end{itemize}
\end{definition}

There is an obvious forgetful map $F\colon\surfo\to\surf$,
that induces a map $F\colon\OA(\surfo)\to\OA(\surf)$.
We choose an initial triangulation $\TT$ of $\surfo$
and let $\T=F(\TT)$ be the induced triangulation of $\surf$.
The quiver with potential $(Q_\TT,W_\TT)$ and the associated Ginzburg dga $\Gamma_\TT$ of $\TT$
are defined to be the ones of $\T$.

As $\surf$ provides a topological model for the cluster category $\hua{C}(\Gamma_\TT)$,
$\surfo$ provides a topological model for the finite dimensional/perfect derived categories
$\D_{fd}(\Gamma_\TT)\subset\per\Gamma_\TT$.
Similar to the case in \S~\ref{sec:ST},
any simple in the canonical heart
$\h_\TT\colon=\mod\Gamma_\TT$ is spherical.
Let
\[
    \ST(\Gamma_\TT)\colon=\< \twi_S \mid S\in\Sim\h_\TT \>
\]
be the spherical twist group
and
\[
    \Sph(\Gamma_\TT)\colon=\ST(\Gamma_\TT)\cdot\Sim\h_\TT
\]
be the set of reachable spherical objects.

\begin{theorem}\cite{QQ,QZ2}\label{thm:QQQZ}
There is a bijection
\begin{gather}\label{eq:X}
\begin{array}{rcl}
    X\colon\CA(\surfo) &\to& \Sph(\Gamma_\TT)/[1],\\
    \eta &\mapsto& X_\eta[\ZZ],
\end{array}
\end{gather}
with the formula
\begin{gather}\label{eq:X int=dim}
    2\Int(\eta_1,\eta_2)=\dim\Hom^\bullet(X_{\eta_1},X_{\eta_2}).
\end{gather}
\end{theorem}

Note that \eqref{eq:X} is compatible with the fact that $\D_{fd}(\Gamma_\TT)$
is CY-3 as interchanging $\eta_i$ does not effect the value on both hand sides.

Moreover, we have the following.
\begin{theorem}\cite{QQ2,QZ2}\label{thm:QQ2QZ2}
There is a bijection
\begin{gather}\label{eq:M2}
\begin{array}{rcl}
    \widetilde{M}\colon\OA^\circ(\surfo) &\to& \Ind^\circ\per\Gamma_\TT,\\
    \widetilde{\gamma} &\mapsto& \widetilde{M}_{\widetilde{\gamma}},
\end{array}
\end{gather}
where $\OA^\circ(\surfo)$ is the set of open arcs
that appear in triangulation obtained from $\TT$ by repeated flipping
and $\Ind^\circ\per\Gamma_\TT$ is a class of rigid indecomposables in $\per\Gamma_\TT$,
with the formula
\begin{gather}\label{eq:M2 int=dim}
    \Int(\widetilde{\gamma},\eta)=\dim\Hom^\bullet(\widetilde{M}_{\widetilde{\gamma}},X_{\eta}).
\end{gather}
for $\widetilde{\gamma}\in\OA^\circ(\surfo)$ and $\eta\in\CA(\surfo)$.
\end{theorem}

\subsection{On the spherical twists: surface case}\label{sec:ST2}
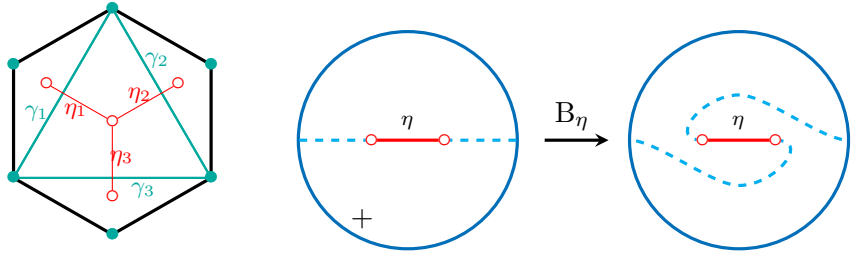
\begin{figure}[t]\centering
\begin{tikzpicture}[scale=.25,rotate=-120]
\foreach \j in {1,...,6}{\draw[very thick](60*\j+30:6)to(60*\j-30:6);}
\foreach \j in {1,...,6}{\draw[Emerald,thick](120*\j-30:6)node{$\bullet$}to(120*\j+90:6)
    (120*\j-90:6)node{$\bullet$};}
\foreach \j in {1,...,3}{  \draw[red](30+120*\j:4)to(0,0);}
\foreach \j in {1,...,3}{  \draw(30-120*\j:4)node[white]{$\bullet$}node[red]{$\circ$};
    \draw(30-120*\j+14:2)node[red]{\text{\footnotesize{$\eta_\j$}}};
    \draw(30-120*\j+24:4)node[Emerald]{\text{\footnotesize{$\gamma_\j$}}};}
\draw(0,0)node[white]{$\bullet$}node[red]{$\circ$};
\end{tikzpicture}
\qquad
\begin{tikzpicture}[scale=.24]
  \draw[very thick,NavyBlue](0,0)circle(6)node[above,black]{$_\eta$};
  \draw(-120:5)node{+};
  \draw(-2,0)edge[red, very thick](2,0)  edge[cyan,very thick, dashed](-6,0);
  \draw(2,0)edge[cyan,very thick,dashed](6,0);
  \draw(-2,0)node[white] {$\bullet$} node[red] {$\circ$};
  \draw(2,0)node[white] {$\bullet$} node[red] {$\circ$};
  \draw(0:7.5)edge[very thick,->,>=stealth](0:11);\draw(0:9)node[above]{$\Bt{\eta}$};
\end{tikzpicture}\;
\begin{tikzpicture}[scale=.24]
  \draw[very thick, NavyBlue](0,0)circle(6)node[above,black]{$_\eta$};
  \draw[red, very thick](-2,0)to(2,0);
  \draw[cyan,very thick, dashed](2,0).. controls +(0:2) and +(0:2) ..(0,-2.5)
    .. controls +(180:1.5) and +(0:1.5) ..(-6,0);
  \draw[cyan,very thick,dashed](-2,0).. controls +(180:2) and +(180:2) ..(0,2.5)
    .. controls +(0:1.5) and +(180:1.5) ..(6,0);
  \draw(-2,0)node[white] {$\bullet$} node[red] {$\circ$};
  \draw(2,0)node[white] {$\bullet$} node[red] {$\circ$};
\end{tikzpicture}

  \caption{The dual graph of a triangulation and the braid twist along a closed arc}
  \label{fig:1}
\end{figure}

The main application of Theorem~\ref{thm:QQQZ} is the following,
which identifies the spherical twist group
with a subgroup, the braid twist group $\BT(\TT)$,
of the mapping class group $\MCG(\surfo)$ of $\surfo$.
More precisely, $$\BT(\TT)\colon=\<B_{\eta}\mid\eta\in\TT^*\>$$
is generated by the braid twists $B_\eta$
along the closed arcs $\eta$ in the dual graph $\TT^*$ of the triangulation $\TT$,
cf. Figure~\ref{fig:1}.

\begin{theorem}\cite[Theorem~1]{QQ}\label{thm:QQ}
The bijection $X$ in \eqref{eq:X} induces an isomorphism
\[\begin{array}{rcl}
    \iota_\TT\colon\BT(\TT)&\cong&\ST(\Gamma_\TT)\\
                    B_\eta&\mapsto&\twi_{X_\eta},
\end{array}\]
sending the generators to the corresponding generators,
which both can be indexed by vertices of $Q_\TT$ (i.e. open arcs in $\TT$).
\end{theorem}

In summary, we have Table~\ref{table} of correspondences.

\begin{table}[h]
\caption{Correspondences II}
\label{table}
\setlength{\extrarowheight}{2pt}
\begin{tabular}{ccc}
\\Topology&&Category\\[1ex]
\hline
Braid twist group $\BT(\TT)$&$\xymatrix@C=2pc{\ar[r]^{\cong}&}$&Spherical twist group $\ST(\Gamma_\TT)$\\
\Huge{$_\curvearrowright$}&&\Huge{$_\curvearrowright$}
\\
Simple closed arcs in $\surfo$&$\xymatrix@C=2pc{\ar[r]^{\text{\tiny{1-1}}}_{\text{\tiny{up to $[1]$}}}&}$
&Spherical obj. in $\D_{fd}(\Gamma_\TT)$\\
\footnotesize{Dual Tri. with Whitehead move} &&
\footnotesize{Hearts with simple tilting}
\\ \hline\hline
\text{\tiny{graph dual}}$\;\;\,\mathrel{\Bigg\updownarrow}\qquad\quad\;\;$
&&$\;\;\;\qquad\qquad\mathrel{\Bigg\updownarrow}\;\;\;$\text{\tiny{sim.-proj. dual}}
\\ \hline\hline
Reachable open arcs in $\surfo$&$\xymatrix@C=2pc{\ar[r]^{\text{\tiny{1-1}}}&}$&
Reachable ind. in $\per\Gamma_\TT$\\
\footnotesize{Triangulations with flip}&&
\footnotesize{Silting objects with mutation}\\
\hline\hline
\text{\tiny{forgetful map}}
$\begin{smallmatrix}\\\surfo\\\mathrel{\bigg\downarrow}\\ \surf\end{smallmatrix}
\qquad\qquad$
&&
\text{\tiny{}}$\qquad\qquad
\begin{smallmatrix}\\\per\Gamma_\TT\\\mathrel{\bigg\downarrow}\\ \C(\Gamma_\TT)\end{smallmatrix}$
\text{\tiny{quotient map}}\\[4ex]
\hline\hline
Open arcs in $\surf$&$\xymatrix@C=2pc{\ar[r]^{\text{\tiny{1-1}}}&}$&
Rigid ind. in $\C(\Gamma_\TT)$\\
\footnotesize{Triangulations with flip}&&
\footnotesize{Cluster tilting objects with mutation}\\
\hline\hline
\end{tabular}
\end{table}

\subsection{Topological Fukaya categories}\label{sec:HKK}
\newcommand{\TFuk}{\operatorname{TFuk}}
There is a closely related result in \cite{HKK},
which is in the same spirit of Theorem~\ref{thm:BZZZZQZ}
and Theorem~\ref{thm:QQQZ}.
Let $S$ be a graded flat surface, $\GA(S)$ the set of graded arcs on $S$,
$\TFuk(S)$ the associated topological Fukaya category
(see \cite[\S~2,\S~3]{HKK} for details).
They prove the following (we only give a simplified version).

\begin{theorem}\cite{HKK}\label{thm:HKK}
There is a bijection
\begin{gather}\label{eq:XX}
\begin{array}{rcl}
    \GA(S) &\to& \Ind\TFuk(S).
\end{array}
\end{gather}
\end{theorem}
The corresponding $\Int=\dim\Hom$ formula will be proved in \cite{IQZ}.

\begin{example}
The prototype for $\TFuk$ is when $S$ is a disk and
$\TFuk(S)\cong\D^b(A_n)$ for an $A_n$ quiver.
\end{example}

\section{Stability conditions via quadratic differentials}\label{sec:Quad}
In this section, we review two important results
realizing spaces of stability conditions as moduli spaces of quadratic differentials.
In both cases, the key idea is the following formula
\begin{gather}\label{eq:idea}
    Z(S)=\int_{\widetilde{\eta}} \sqrt{\phi},
\end{gather}
that the central charge $Z$ of a semistable object $S$ should be given by
the the period (integration) of some differential (i.e. $\sqrt{\phi}$) along some
curve $\widetilde{\eta}$ --- the double cover of a closed arc $\eta$ corresponds to $S$.
The correspondence $$\eta \longleftrightarrow S$$
is in fact the ones in Theorem~\ref{thm:QQQZ} and Theorem~\ref{thm:HKK}.

\subsection{Calabi-Yau-3 case}\label{sec:BS}
\newcommand{\Pol}{\operatorname{Pol}}
\newcommand{\Squad}{\Quad^\times}

Let $\mathbf{X}$ be a compact Riemann surface and $\omega_\mathbf{X}$
its holomorphic cotangent bundle.
A \emph{meromorphic quadratic differential} $\phi$ on $\mathbf{X}$ is a meromorphic section
of the line bundle $\omega_{\mathbf{X}}^{2}$.
We will only consider \emph{signed GMN differentials} $\Phi=(\phi,\zeta)$ on $\mathbf{X}$,
consists of a meromorphic quadratic differential $\phi$ whose zeros are all simple (order is one)
and a function $\zeta\colon\Pol_2(\phi)\to\{\pm1\}$ on the set of order 2 poles of $\phi$.
Denote by $\mathbf{X}^{\Phi}$ the real blow-up of $\mathbf{X}$ at the poles with order at least 3.
Denote by $\Squad(\surf)$ the moduli space of signed quadratic differentials on $\surf$,
consisting of a signed GMN differentials $\Phi$ on some Riemann surface $\mathbf{X}$
and a diffeomorphism $\surf\to\mathbf{X}^\Phi$
(sending $\mathbf{P}$ to the set of poles with order one or two).

In \cite{BS}, they prove the following.
\begin{theorem}\cite[Thm.~11.2]{BS}\label{thm:BS}
As complex manifolds, there is an isomorphism
\begin{equation}\label{eq:BS}
    \Stap(\D_{fd}(\Gamma_\T))/\Aut^\circ\D_{fd}(\Gamma_\T)\cong\Squad(\surf).
\end{equation}
\end{theorem}

When $\mathbf{P}=\emptyset$,
there is no sign involved (so $\Squad(\surf)=\Quad(\surf)$).
Moreover, there are couple of variations (framed version) of \eqref{eq:BS} given in \cite{KQ2}.
In particular, we have the following.

\begin{theorem}\cite[Thm.~4.14 and Thm.~4.16]{KQ2}
Let $\surfo$ be a decorated marked surface with a triangulation $\TT$.
Then
$$\Stap(\Gamma_\TT)\cong\operatorname{FQuad}^\TT(\surfo)$$
is simply connected,
where $\operatorname{FQuad}$ is the $\surfo$-framed version of $\Quad(\surf)$.
\end{theorem}

Also, \cite{I} generalizes Theorem~\ref{thm:BS} to CY-$N$ case when $\surf$ is a polygon
(and the quiver with potential is of type $A_n$).
\subsection{Non-Calabi-Yau case}\label{sec:HKK2}

Building on Theorem~\ref{thm:HKK}, there is a similar result of Theorem~\ref{thm:BS}
as follows.

In \cite{HKK}, they prove the following.
\begin{theorem}\cite[Thm.~5.3]{HKK}\label{thm:HKK2}
As complex manifolds, there is an isomorphism
\begin{equation}\label{eq:q=s}
    \Quad(S)\cong \Stab^\bullet(\TFuk(S)),
\end{equation}
where $\Quad(S)$ is the $S$-framed quadratic differential with exponential type singularities
and $\Stab^\bullet(\TFuk(S))$ is the union of certain connected components of $\Stab(\TFuk(S))$.
\end{theorem}

Note that when $S$ is a disk with three marked points/boundary components,
$$\TFuk(S)\cong\D^b(A_2)$$ and the result above is shown in \cite{BQS}.
\subsection{Further studies}\label{sec:IQ}
Clearly, the results in the two cases, namely,
\begin{itemize}
\item Theorem~\ref{thm:QQQZ} for category and Theorem~\ref{thm:BS} for geometry in CY-3 case;
\item Theorem~\ref{thm:HKK} for category and Theorem~\ref{thm:HKK2} for geometry in non-CY case,
\end{itemize}
are of the same sort.
In the upcoming series of works \cite{IQ1,IQ2,IQZ},
we will introduce new categories and $q$-deformation of stability conditions
to relate them in precise statements, both categorically and geometrically.


\end{document}